# Union-closed Sets Conjecture Holds for Height H($\mathcal{F}$)≤ 3 and H($\mathcal{F}$)≥ n − 1

Chenxiao Tian[1]

*Department of mathematics, Peking University, No. 5 Summer Palace Road, Haidian District, Beijing, China*

**Abstract**

For each given union-closed family $\mathcal{F}$(n, m) of n elements and m sets, we discuss the union-closed sets conjecture from height number of the UC family, which is a natural parameter from lattice theory. In this paper, we call it height number of $\mathcal{F}$(n, m), recorded as H($\mathcal{F}$). we prove that for any given union-closed family $\mathcal{F}$, union-closed sets conjecture holds if its height number H($\mathcal{F}$)≤ 3 or H($\mathcal{F}$)≥ n − 1.

Since the height number H($\mathcal{F}$) is a positive integer which is bounded between 1 to n. As an attempted approach and framework, if we can prove union-closed sets conjecture holds for all possible value of H($\mathcal{F}$), then union-closed sets conjecture is true.

————————————————————————————————————————

## 1. Introduction

Union-closed family $\mathcal{F}$ (n, m) is a finite family of n elements and m sets that is closed under the union operation of sets. In this paper, we will define a global character number H($\mathcal{F}$) of a given union-closed family, we name it height number in this paper and will show how it may lead an attempted approach to union-closed sets conjecture, which is a celebrated conjecture stated by Frankl in 1979:

**Union-closed sets conjecture.** For any finite union-closed family of sets $\mathcal{F}$ without ∅ there exists an element x which belongs to more than half of the sets of $\mathcal{F}$.

In this paper, we continuously assume that for any given finite union-closed family of sets $\mathcal{F}$, it doesn't contain an empty set {∅} without special statement. If we allow it to contain an empty set, then the conjecture's assertion needs to be changed into "there exists an element x which belongs to at least half of the sets of $\mathcal{F}$".

The structure of the paper is the following, in the first part of the paper we will define what the height number of a union-closed family $\mathcal{F}$(n, m) of n elements and m sets is, we record it as H($\mathcal{F}$). Then in the second part of the paper, we will prove that union-closed sets conjecture holds for union-closed family $\mathcal{F}$ whose height number H($\mathcal{F}$)≤3 or H($\mathcal{F}$)≥ n − 1. In this part, it will be divided into four cases, the first case is for H($\mathcal{F}$)=1 or H($\mathcal{F}$)=2. The second case is for H($\mathcal{F}$)=3 and m>2n. The third case is for H($\mathcal{F}$)=3 and m≤2n. The fourth case is for H($\mathcal{F}$) ≥ n − 1.This main theorem of this paper leads to an establishment of the union-closed sets conjecture for a relatively large class of UC-family. And since H($\mathcal{F}$) is a positive integer which is bounded between 1 to n. It may lead a potential way to approach union-closed sets conjecture in the future work.

[1] Email: 1700013239@pku.edu.cn

## 2. The Height Number of Union-closed Family

For completeness and some extra terms, we will define the height from lattice theory for union-closed family specifically. Firstly, we describe a decomposition operation of a given union-closed family $\mathcal{F}$:

In any given union-closed family $\mathcal{F}$ of sets, we define the subfamily $\pi_1$ of $\mathcal{F}$ to be the following connection of sets in $\mathcal{F}$:

$\pi_1 = \{A \in \mathcal{F} | \forall B \in \mathcal{F}, B \neq A \to B \not\subset A\}$

If $\mathcal{F}$ still has other sets which not contain in $\pi_1$, then we consider the rest sets of $\mathcal{F}$, we record it as $\mathcal{F} \backslash \pi_1$. Similarly, we define $\pi_2$, to be the following connection of sets in $\mathcal{F} \backslash \pi_1$:

$\pi_2 = \{A \in \mathcal{F} \backslash \pi_1 | \forall B \in \mathcal{F} \backslash \pi_1, B \neq A \to B \not\subset A\}$

Continuously, if $\mathcal{F} \backslash (\pi_1 \cup \pi_2)$ is still not empty, we can define $\pi_3$:

$\pi_3 = \{A \in \mathcal{F} \backslash (\pi_1 \cup \pi_2) | \forall B \in \mathcal{F} \backslash (\pi_1 \cup \pi_2), B \neq A \to B \not\subset A\} \backslash$

Since $\mathcal{F}$ is a finite family of sets, we can continuously carry out this process until $\mathcal{F} \backslash (\cup_{i=1}^{H} \pi_i)$ becomes an empty set, where H exactly records how many times we execute this process. Now that $\mathcal{F}$ is a family of finite sets, so we can assure that we can end this decomposition in finite steps.

### Definition 2.1 (Height Decomposition and Height Number)

For any given union-closed family $\mathcal{F}$, we can continuously operate the above decomposition operation and end it in finite steps. We record the number of steps we execute as H($\mathcal{F}$), then:

(i) We define the decomposition $\mathcal{F} = \cup_{i=1}^{H(\mathcal{F})} \pi_i$ as the height decomposition of $\mathcal{F}$.

(ii) We define the positive integer number H($\mathcal{F}$) as the height number of $\mathcal{F}$.

(iii) We name $\pi_{H(\mathcal{F})}$ as the $1^{st}$ height of the union-closed family $\mathcal{F}$, for any $\pi_i$ we name it as the $(H(\mathcal{F}) + 1 - i)^{th}$ height of the union-closed family $\mathcal{F}$.

(iv) We say the $k^{th}$ height is higher than the $k'^{th}$ height if and only if $k < k'$.

(v) We say the height number of a set A in $\mathcal{F}$, that means the order number of the height which A belongs to.

### Remark 2.1

It is obviously that in a union-closed family $\mathcal{F}$ of n elements, we may just assume that $\cup_{A \in \mathcal{F}} A = \{1,2,3..., n\}$, the first height $\pi_{H(\mathcal{F})}$ must be composed of the set $\{1,2,3..., n\}$, that is, $\pi_{H(\mathcal{F})} = \{\{1,2,3 ..., n\}\}$. And if a union-closed family $\mathcal{F}$ is allowed to contain set $\emptyset$ and it does have an empty set, then the $(H(\mathcal{F}))^{th}$ height $\pi_1$ of the union-closed family is $\{\emptyset\}$.

Another obvious fact is that a set B's height number is larger than A is equal to a set A's height is higher than B.

**Property 2.1**

The following six basic properties of height number and heights of a union-closed family is related to the later discussion in section 3:

(i) For any given $k^{th}$ height of a union-closed family $\mathcal{F}(n, m)$, there must exist at least one set which has at most (n+1-k) elements.

(ii) For any given union-closed family $\mathcal{F}(n, m)$, its height number is an positive integer bounded between 1 to n.

(iii) There exists at least one element $x \in \bigcup_{A \in \mathcal{F}} A$, for any given height of this union-closed family, it appears in at least one set in this height.

(iv) If B and A are located in the same height and B≠ A. Then $B \not\subset A$ and $A \not\subset B$.

(v) If $B \subset A$ and B≠ A, then A is in a strictly higher height than B, that is, A's height number is strictly lower than B's height number.

(vi) For any given set in the $2^{th}$ height of a union-closed family $\mathcal{F}(n, m)$, it has at most (n-1) elements.

*Proof:*

(i) We keep assuming that $\bigcup_{A \in \mathcal{F}} A$ ={1,2,3..., n}. For any given UC family, we can always pick at least one chain B which starts from the set C={1,2..., n} located in $\pi_{H(\mathcal{F})}$ and ends in a set D located in $\pi_1$. Now we consider any two sets in this chain which are located in two adjacent heights, let them be X∈ $\pi_i$ and Y∈ $\pi_{i+1}$, then $X \subseteq Y$ and X ≠ Y. So |X| is strictly lower than |Y|. Since in the $k^{th}$ height, the set E in chain B has at most (n+1-k) elements.

(ii) As for property (ii), since we assume that there is no empty set belonging to a given union-closed family in this paper and the $1^{th}$ height $\pi_{H(\mathcal{F})} = \{\{1,2,3 \dots, n\}\}$, it is a direct corollary of property (i).

(iii) Since we assume that ∅ doesn't belong to the union-closed family $\mathcal{F}$, (iii) is a direct corollary from the proof of (i).

(iv) Otherwise, if B⊂A or A⊂B, then A and B must be located in different heights, which makes a contradiction.

(v) According to the property (iv), they are not in the same height. And we look back at the construction of height, A will be classified in a height later than B, in other words, its height number is lower than B and its height is higher than B.

(vi) It's obvious that except {1,2..., n}, other sets' size in family are all no more than (n-1).

## 3. The Union-closed sets conjecture holds for H($\mathcal{F}$)≤ 3 and H($\mathcal{F}$)≥ $n-1$

Now we go to the main theorem of this paper. The proof will be divided into three parts, firstly we will study the case for union-closed family $\mathcal{F}$ (n, m) whose H($\mathcal{F}$)=1or H($\mathcal{F}$)=2, particularly we will study the case for H($\mathcal{F}$)=2 because we will see that union-closed family whose height number is 2 plays an important role in the construction of general union-closed family later. Then we will study the case for H($\mathcal{F}$)= 3, it will also be divided into another two cases. The first case is that H($\mathcal{F}$)= 3 and m >2n, the second case is that H($\mathcal{F}$)= 3 and m ≤2n. Finally, we prove that it also holds for H($\mathcal{F}$)≥ n − 1.

Before we start our proof, we first state the main theorem of this section as following:

**Theorem 3.1**

For any finite union-closed family of sets $\mathcal{F}$ (n, m), if its height number H($\mathcal{F}$) ≤ 3 or H($\mathcal{F}$)≥ n − 1, then there exists an element x which belongs to more than half of the sets of $\mathcal{F}$.

## 3.1 The Case For H($\mathcal{F}$)=1 and H($\mathcal{F}$)=2

Firstly, the case for **H($\mathcal{F}$)=1** is trivial, we list it as a simple fact for the completeness of the proof.

**Fact 3.1.1**

For any finite union-closed family of sets $\mathcal{F}$ (n, m), if its H($\mathcal{F}$)= 1, then m=1 and it is simply the form of {{1,2,3…, n}}. So every element in $\mathcal{F}$ belongs to more than half of the sets of $\mathcal{F}$.

Then we come to the case for H($\mathcal{F}$)=2, we will prove the following lemma and union-closed sets conjecture for H($\mathcal{F}$)=2 is a direct corollary of this lemma

**Lemma 3.1.1**

For any finite union-closed family of sets $\mathcal{F}$ (n, m), if its H($\mathcal{F}$)=2, then it has following two properties:

(i) For any given element x ∈ $\bigcup_{A \in \mathcal{F}}$ A. The number of sets in the $2^{nd}$ height $\pi_1$ that doesn't contain element x is no more than 1.
(ii) The number of sets in the $2^{nd}$ height $\pi_1$ is bounded between 1 to n.

*Proof:*
(i) Otherwise, if there exist two sets A and B that doesn't contain element x. We consider the union of A and B, let C=A ∪ B. Since it is a union-closed family, so C is also in $\mathcal{F}$. However, C doesn't contain element x since both of A and B don't contain element x. So the set C doesn't belong to the first height of $\mathcal{F}$, that is $\pi_2$ = {{1,2, … , n}}. And on the other hand, since A ⊂ C and B ⊂ C, C doesn't belong to the second height of $\mathcal{F}$. Now that C doesn't belong to the first height and the second height, so C doesn't belong to $\mathcal{F}$. This makes a contradiction.

(ii) Otherwise, if there exist more than n sets in the $2^{nd}$ height. According to the lemma 3.1.1 (i), since each element can at most not appear in more than one set in the second height, then if there are more than n sets in the second height, there must exist a set that contain all n elements of this union-closed family. That is, there exists a set {1,2 … , n} in the second height in the union-closed family. This also makes a contradiction.

From lemma 3.1.1, we directly get the following corollary:

**Corollary 3.1.1**

For any finite union-closed family of sets $\mathcal{F}$ (n, m), if its H($\mathcal{F}$)=2, then any element in this union-closed family appears at least in more than m-1 sets. Since we can select an element x which appears in more than half of the sets in this family.

*Proof:*
If the second height of this union-closed family contains only one set, according to property 2.1 (iii), we can ensure that we can select an element x which appears in both of the sets in this family, which is more than half sets of the family.

If the second height of this union-closed family contains more than one set, then every element in this union-closed family appears in at least m-1 sets of this family, since m>2, so we get $\frac{m-1}{m} > \frac{1}{2}$. It implies that it appears in more than half of the sets in this family.

## 3.2 The Case For H($\mathcal{F}$)=3 and $m > 2n$

In this subsection, we are going to prove that union-closed sets conjecture is true for the case H($\mathcal{F}$)=3 and $m > 2n$. Before we start our proof, we firstly are going to focus on some important subfamily of a union-closed family $\mathcal{F}$ which are also union-closed and their height number H($\mathcal{F}$)=2. We will see that this kind of union-closed family plays an important role in constructing general union-closed family, for more concise and imagery description in the later proof, we temporarily use tent to represent this sub structure.

**Definition 3.2.1** (**Tent Sub Union-closed Family**)
Given a finite union-closed family of sets $\mathcal{F}$ (n, m) whose height number H($\mathcal{F}$)≥2, let set A and B both be sets of the UC family which don't belong to the (H($\mathcal{F}$))$^{th}$ height of the family and A is located in the k$^{th}$ height. Then a tent sub union-closed family is defined as following:
  (i) We say T is a tent sub union-closed family if and only if T is a union-closed family and its H($\mathcal{F}$)=2. We call it short for tent in this paper.
  (ii) We say T(A) is a tent induced by set A, that means T(A) includes A itself and all the sets in the $(k+1)^{th}$ height which are included in A.
  (iii) If A and B are both in the k$^{th}$ height and $A \neq B$ ( The height should have at least two sets for the reasonability of this definition.), then:
    We say the intersection number Int (A, B) of T(A) and T(B) is the number of sets in the $(k+1)^{th}$ height which both belong to T(A) and T(B), that is, Int (A, B) = |T(A) ∩ T(B)|.

The following lemma explains the reasonability of definition 3.2.1 (ii).

**Lemma 3.2.1**
T(A) defined as above is also a union-closed family whose H($\mathcal{F}$)=2.

*Proof:*
If T(A) contains only one set in the $(k+1)^{th}$ height, the claim is already true.
Else if T(A) contains more than one set, we select two of them arbitrarily, let them be B and C. We consider their union $B \cup C$. If we assume that $B \cup C$ is not equal to A, according to property 2.1 (iv), $B \cup C$ can't be located in the $(k+1)^{th}$ height or the k$^{th}$ height of the union-closed family. And according to property 2.1 (v), now that B ∪ C is strictly included in A, so its height number must be strictly higher than k. We also notice that B is strictly included in B ∪ C for property 2.1 (iv), so its height number must be strictly lower than k+1. This makes a contradiction.

So from lemma 3.2.1, we can see the reasonability of definition 3.2.1 (ii).

### Corollary 3.2.1

Since tent T is also a union-closed family whose height number H($\mathcal{F}$)=2, so it keeps all the property of union-closed family whose H($\mathcal{F}$)=2, which implies that for each element in $\bigcup_{A \in T} A$, it appears in at least $(|T| - 1)$ sets of this tent (Where |X| represents the number of elements in the set X).

The following lemma gives a bound for the intersection number between two different tents induced by two sets in a same height.

### Lemma 3.2.2

Let A and B be different sets of a given union-closed family in the same height which is higher than $\pi_1$. Of course, the height has at least 2 sets for well-defined purpose. Then the intersection number Int (A, B) = |T(A)∩T(B)| is on more than 1, that is, Int (A, B)≤ 1.

*Proof:*

If we assume that T(A)∩T(B) has two different sets, let them be C and D. Then we consider C ∪ D. Since C and D both belong to T(A) and T(A) is a tent according to lemma 3.2.1, so we have C ∪ D=A. Similarly, we have C ∪ D=B when we consider tent T(B). Since A=B, this makes a contradiction.

Now we are going to prove the main conclusion of this section.

### Proposition 3.2.1

For any finite union-closed family of sets $\mathcal{F}$ (n, m), if the height number H($\mathcal{F}$) ≤ 3 and m > 2n. Then there must exist an element x which belongs to more than half of the sets of $\mathcal{F}$.

*Proof:*

In this case, we select the element $x \in \bigcup_{A \in \mathcal{F}} A$ to be the element which appears $(|\pi_2| - 1)$ times in the second height. This element does exist, because the union $\pi_3 \cup \pi_2$ of the first height and the second height is also a tent, according to Corollary 3.1.1, the element x appears in at least $(|\pi_2| - 1)$ times in the second height and $|\pi_2|$ times in the whole tent. However, there exist at least one element x that appears $(|\pi_2| - 1)$ times and only $(|\pi_2| - 1)$ times in the second height $\pi_2$, otherwise all the sets in $\pi_2$ must be {1, 2, ..., n}, of course this is a contradiction.

Now let M be the set in $\pi_2$ which doesn't contain the element x, let T(M) be the tent induced by the set M. We will prove that all the sets in $\pi_1 \setminus T(M)$ must have the element x.

Firstly, if $\pi_1 \setminus T(M)=\emptyset$, then the whole family has at most 1+n+n-1=2n sets, this makes a contradiction to $m > 2n$.

Then for any given set $A \in \pi_1$ and $A \notin T(M)$, we can pick a set B in $T(M) \cap \pi_1$. It is obvious that $B \in T(M)$ doesn't contain the element x. We consider $A \cup B$. Then if A also doesn't contain the element x, then $A \cup B$ doesn't contain the element x. Now that according to property2.1(iv) and (v), $A \cup B$ must be located in $\pi_2$ or $\pi_3$. The only set in $\pi_2$ and $\pi_3$ that doesn't have the element x is exactly M, so we get $A \cup B = M$. Since we have $A \subset M$, that is, $A \in T(M)$. This makes a contradiction to $A \notin T(M)$. So we get that all the sets in $\pi_1 \setminus T(M)$ must contain the element x.

So we finally get that the number of sets that don't contain the element x in this whole union-closed family $\mathcal{F}$ is exactly $|T(M)|$. According to lemma 3.1.1(ii) and property 2.1(vi) again, $|T(M)| \leq n$. Now we assume that $m > 2n$, so we get that $\frac{m-|T(M)|}{m} = 1 - \frac{|T(M)|}{m} \geq 1 - \frac{n}{m} > \frac{1}{2}$. This completes the whole proof of proposition 3.2.1.

## 3.3 The Case For H($\mathcal{F}$)=3 and $m \leq 2n$

In this section, for the completeness we firstly recall the definition of a separating union-closed family.

### Definition 3.3.1 (Separating Union-closed Family)
A finite union-closed family $\mathcal{F}$ is a separating union-closed family if and only if for each pair of different elements x and y, there always exists at least one set A which belongs to $\mathcal{F}$ that doesn't simultaneously contain these two different elements.

It's obvious that we have the following fact, which is an equal form of union-closed sets conjecture.

### Fact 3.3.1
The Union-closed sets conjecture holds for height no more than k if and only if it holds for all separating Union-closed families whose height no more than k.

About separating union-closed family we have the following result which is mentioned in Henning Bruhn and Oliver Schaudt's literature review ([BS15]), which is marked as theorem 23.By their introduction, it is an analogous result for the small family compared to the similar result for the very large family obtained by Nishimura and Takahashi([TS96]). And it can be directly proved by the construction used in the Falgas-Ravry's theorem ([FR11]) aiming at estimating the average frequency of elements in a separating union-closed family $\mathcal{F}$.

### Theorem 3.3.1 (Theorem 23 in [BS15])
Any separating family on n elements with at most 2n member-sets satisfies the union-closed sets conjecture.

### Corollary 3.3.1

For any finite union-closed family of sets $\mathcal{F}$ (n, m), if the height number H($\mathcal{F}$) ≤ 3. Then there must exist an element x which belongs to more than half of the sets of $\mathcal{F}$.

*Proof:*
Firstly, we can get a direct corollary from proposition 3.2.1, that is, the Union-closed sets conjecture holds for any separating finite union-closed family of sets $\mathcal{F}$ (n, m) whose H($\mathcal{F}$)≤ 3 and m＞2n . According to fact 3.3.1 and the conclusion of theorem 3.3.1 , we put them together and prove the corollary 3.3.1.

### 3.4 The Case For H($\mathcal{F}$)≥ n − 1

Now we prove that it holds for H($\mathcal{F}$)≥ n − 1, firstly we recall the definition of FC-family:
**Definition 3.4.1(FC-family (see [Va02] or [Po92]))**
An FC-family is a UC family B such that for every UC-family A, if B⊆A then A satisfies the Union-closed sets conjecture.
**Proposition 3.4.1**
For any finite union-closed family of sets $\mathcal{F}$ (n, m), if its height number H($\mathcal{F}$)≥ n − 1. Then there must exist an element x which belongs to more than half of the sets of $\mathcal{F}$.
*Proof:*
If H($\mathcal{F}$)=n, according to property 2.1 (i), there must exist at least one set that contain only one element in $\pi_1$, let it be {x}, this is a known FC-family

If H($\mathcal{F}$)=n-1, according to property 2.1 (i), there must exist at least one set that contain only one or two elements in $\pi_1$, let it be {x} or {x,y}, this is also a known FC-family.

To put it in a nutshell, from corollary 3.3.1, proposition 3.2.1, corollary 3.1.1, fact 3.1.1, proposition 3.4.1. We finally finish the whole proof of theorem 3.1.

### Acknowledgement

I thank Bruhn, H. for his important suggestions on my paper, which helps me to give a more concise writing on this paper. I also thank Bruhn, H and Schaudt, O for their fantastic literature review, which helps me a lot in this paper's writing and gives me an overview of this interesting conjecture. I also thank my undergraduate school Peking University, my research is supported by its top undergraduate student research training subsidy.